\documentclass{ifacconf}
\usepackage{graphicx}      
\usepackage{natbib}        
\usepackage{color}

\usepackage[hang,small]{caption}
\usepackage{amsmath}
\usepackage{amssymb}
\usepackage{amsfonts}

\usepackage{algorithm}
\usepackage{algpseudocode}

\DeclareMathOperator*{\argmin}{arg\,min}

\newcommand{\bx}{\mathbf{x}}
\newcommand{\pd}{\partial}

\newcommand{\CC}{\mathcal{C}}
\newcommand{\CN}{\mathcal{N}}
\newcommand{\CD}{\mathcal{D}}
\newcommand{\CE}{\mathcal{E}}

\newcommand{\BRnd}{\mathbb{R}^{nd}}

\newcommand{\BRCd}{\mathbb{R}^{|\mathcal{C}_l|d}}
\newcommand{\clqG}{\mathrm{clq}(G)}
\newcommand{\clqiG}{\mathrm{clq}_i(G)}

\newcommand{\CX}{\mathcal{X}}
\newcommand{\hl}{\hat{l}}

\begin{document}
\begin{frontmatter}

\title{Accelerated Distributed Projected Gradient Descent for Convex Optimization with Clique-wise Coupled Constraints
} 

\thanks[footnoteinfo]{This work was partially supported by the joint project of Kyoto University and Toyota Motor Corporation, titled “Advanced Mathematical Science for Mobility Society”.}

\author[First]{Yuto Watanabe} 
\author[Second]{Kazunori Sakurama} 

\address[First]{Department of System Science, Graduate school of Informatics, Kyoto University (email: y-watanabe@sys.i.kyoto-u.ac.jp)}
\address[Second]{Department of System Science, Graduate school of Informatics, Kyoto University(email: sakurama@sys.i.kyoto-u.ac.jp)}

\begin{abstract}    
This paper addresses a distributed convex optimization problem with a class of coupled constraints, which arise in a multi-agent system composed of multiple communities modeled by \textit{cliques}.
First, we propose a fully distributed gradient-based algorithm with a novel operator inspired by the convex projection, called the clique-based projection.
Next, we scrutinize the convergence properties for both diminishing and fixed step sizes.
For diminishing ones, we show the convergence to an optimal solution under the assumptions of the smoothness of an objective function and the compactness of the constraint set.
Additionally, when the objective function is strongly monotone, the strict convergence to the unique solution is proved without the assumption of compactness.
For fixed step sizes, we prove the non-ergodic convergence rate of $O(1/k)$ concerning the objective residual under the assumption of the smoothness of the objective function.
Furthermore, we apply Nesterov's acceleration method to the proposed algorithm and establish the convergence rate of $O(1/k^2)$.
Numerical experiments illustrate the effectiveness of the proposed method.
\end{abstract}


\end{frontmatter}

\section{Introduction}

Due to the rapid improvement of information technology, devices with built-in computers are becoming more prevalent in various application domains (e.g., sensor networks, traffic networks, and power systems).
As a result, systems are getting larger and larger by connecting these devices.
In such large-scale systems, many challenges arise.
For example, the central management of these systems leads to heavy computational burdens, and low-spec devices are incorporated into the network.
To deal with these challenges, researchers in the control and machine learning communities have vigorously investigated distributed optimization, which aims to solve optimization problems by local communication without central management.
Various papers have addressed the issue of the fast convergence and efficient computation of algorithms to reduce computational costs, communication frequency, and so on.

In recent years, distributed optimization problems with constraint-coupling have been studied in several papers (e.g., \cite{Falsone2020-mr,Notarnicola2020-bn,Chang2016-lz,Su2022-tn,Wu2022-vn}).
As opposed to conventional methods like the dual decomposition in \cite{Terelius2011-ir} and the ADMM in \cite{Boyd2011-yu}, the methods proposed in those papers do not require central management.
Such a constraint-coupled problem is a practical and general framework, which contains 
the widely studied problem in which the agents have to seek a common solution under agent-wise constraints, e.g., \cite{Nedic2010-xf,Zhu2012-am,Yang2019-uc}, etc.
However, the existing methods for constraint-coupling setups have slow convergence rates or do not explicitly show their convergence rate.
\cite{Falsone2020-mr} proposed an ADMM-based algorithm with the dynamic average consensus strategy but did not explicitly show the convergence rate.
Similarly, \cite{Notarnicola2020-bn} proposed a distributed algorithm based on a relaxed primal problem and the duality theory without offering the convergence rate.
In \cite{Chang2016-lz,Wu2022-vn}, the convergence rate of $O(1/k)$ was established in an ergodic sense.
Just recently, in \cite{Su2022-tn}, the non-ergodic convergence rate of $O(1/k)$ has been achieved for constraint-coupled optimization problems.
To the authors' knowledge, the convergence rate of existing methods for constraint-coupling setups has been at most $O(1/k)$ so far.
Given implementation to low-spec devices, it is meaningful to construct a faster and more efficient method.

This paper addresses a distributed convex optimization problem over a multi-agent network with a class of coupling in constraints called clique-wise coupled constraints.
This class corresponds to a system consisting of multiple communities (with some overlap) modeled by cliques, i.e., complete subgraphs (see \cite{Bollobas1998-gp}), and each community has some constraints.
This formulation contains many conventional distributed optimization problems, e.g., the consensus and agent-wise constraints studied in \cite{Nedic2010-xf}, etc.
First, we develop a clique-based projection operator, which is a novel extension of the convex projection.
Second, by using this developed operator, we propose a distributed optimization algorithm,
the \textit{clique-based projected gradient descent} (CPGD), which is a generalization of the well-known projected gradient descent (PGD, see \cite{Calamai1987-ai}).
Next, we prove its convergence to an optimal solution for diminishing and  fixed step sizes.
Moreover, using the Nesterov acceleration scheme (\cite{Nesterov1983-ir,Beck2009-bq}), we show that the proposed CPGD achieves the outstanding convergence rate of $O(1/k^2)$.
Finally, we demonstrate the effectiveness of the proposed method through numerical experiments.

The major contributions of our proposed method, CPGD are threefold as follows.
(i) The proposed CPGD does not require central management and is implementable only with peer-to-peer communication between agents.
(ii) In the CPGD, by repeatedly operating the clique-based projection, we can generate a sequence arbitrarily close to a constraint set, which is challenging for the existing methods such as \cite{Falsone2020-mr,Notarnicola2020-bn}, etc.
Although there is a trade-off between the number of its operation and communication costs, the number can be tuned depending on the situation.
(iii) The CPGD is fast and efficient for some practical constraints.
The CPGD with the Nesterov acceleration achieves the non-ergodic convergence rate $O(1/k^2)$, which is faster than the existing distributed algorithms for constraint-coupled setups.
Moreover, the computational cost in each iteration is very small for some typical constraints (e.g., linear and norm constraints).

The rest of this paper is organized as follows.
Section \ref{sec:pre} provides preliminaries.
Section \ref{sec:problem} presents the problem setting and show several examples that this paper considers.
In Section \ref{sec:algorithm}, we extend the convex projection and propose a new distributed algorithm, the clique-based projected gradient descent (CPGD).
In Section \ref{sec:convergence}, we show the convergence properties of the CPGD.
Section \ref{sec:acc} presents the accelerated CPGD, which accomplishes the convergence rate of $O(1/k^2)$.
In Section \ref{sec:numerical_experiments}, numerical examples illustrate the effectiveness of our CPGD.
Finally, Section \ref{sec:conclusion} concludes this paper.

\section{
Preliminaries
}\label{sec:pre}
\subsection{Notation}\label{subsec:notation}
Let $\mathbb{R}$ and $\mathbb{N}$ be the set of real numbers and that of positive integers, respectively.
Let $|\cdot|$ be the number of elements in a countable finite set.
The closure of a set is denoted by $\mathrm{cl}(\cdot)$.
For a mapping $T:\mathbb{R}^m\to\mathbb{R}^m$, define the fixed points set of $T$ as $\mathrm{Fix}(T)=\{x\in\mathbb{R}^m:T(x)=x\}$.
Let $I_d\in \mathbb{R}^{d\times d}$ denote the $d\times d$ identity matrix.
Let $\mathbf{1}_d=[1,\ldots,1]^\top\in \mathbb{R}^d$ denote the vector of $d$ ones.
With a positive definite and symmetric matrix $Q\in\mathbb{R}^{m\times m}$, we define the norm $\|\cdot\|_Q$ as $\| v \|_Q=\sqrt{v^\top Q v}$ for a vector $v\in\mathbb{R}^m$.
When $Q=I_d$, we simply write $\|\cdot\|_{I_d}$ as $\|\cdot\|$.

For a vector $v=[v_1^\top,\ldots,v_j^\top,\ldots,v_{N}^\top]^\top\in\mathbb{R}^{N d}$ with vectors $v_1,\ldots,v_{N}\in\mathbb{R}^d$, $[v]_j$ represents the operation to extract the $j$th vector $v_j$ from $v$, that is,
\begin{align*}
   [v]_j=v_j\in\mathbb{R}^{d}.
\end{align*}

For a vector $x=[x_1^\top,\ldots,x_n^\top]^\top\in\BRnd$ with $x_1,\ldots,x_n\in\mathbb{R}^d$ and  a subset $\CC=\{j_1,\ldots,j_{|\CC|}\}\subset\{1,\ldots,n\}$, let $x_{\CC}$ be 
\begin{equation*}
    x_\CC = [x_{j_1}^\top,\ldots,x_{j_{|\CC|}}^\top]^\top \in \mathbb{R}^{|\CC|d},
\end{equation*}
where $\{j_1,\ldots,j_{|\CC|}\}$ is a strictly monotonically increasing sequence.

For $x=[x_1^\top,\ldots,x_n^\top]^\top\in\BRnd$ and a differentiable function $f:\BRnd\to\mathbb{R}$, we write $\nabla f(x)=[\nabla_1 f(x)^\top,\ldots,\nabla_n f(x)^\top]^\top \in \BRnd$ with $\nabla = \pd /\pd x$ and $\nabla_i=\pd/\pd x_i.$

\subsection{Graph Theory}
In this subsection, we provide graph theoretic concepts.
Consider a graph $G=(\CN,\CE)$ with a node set $\CN=\{1,\ldots,n\}$ and an edge set $\CE$ consisting of pairs $(i,j)$ of nodes $i,j\in \CN$.
If $(i,j)\in\CE \Leftrightarrow (j,i)\in\CE$ holds for all $(i,j)\in\CE$, the graph $G$ is said to be \textit{undirected}.
In the following, we consider a time-invariant undirected graph $G$.
For $i\in\CN$ and $G$, let $\CN_i\subset \CN$ be the \textit{neighbor set} of node $i$ over $G$, defined as $\CN_i=\{j\in \mathcal{N}:(i,j)\in \mathcal{E}\}\cup\{i\}$.

For an undirected graph $G$,
we consider a set $\CC\subset\CN$.
For $\CC$ and $\CE$,  let $\CE|_{\CC}$ denote the subset of $\CE$ defined as $\mathcal{E}|_{\mathcal{C}}=\{(i,j)\in \mathcal{E}:i,j\in \mathcal{C}\}$. We call $G|_\CC=(\CC,\CE|_\CC)$ a subgraph induced by $\CC$.
If $G|_\CC$ is complete, $\CC$ is called a \textit{clique} in $G$.
If a clique $\CC$ is not contained by any other cliques, $\CC$ is said to be \textit{maximal}.
Let $\mathrm{clq}(G)=\{1,2,\ldots,q\}$ be a set of indices of maximal cliques in $G$. 
For $i\in\CN$, we define $\mathrm{clq}_i(G)$ as an index set of the maximal cliques containing $i$, that is, $\mathrm{clq}_i(G)=\{k\in \mathrm{clq}(G):i\in \mathcal{C}_k \}$.
Note that, for each $i\in\CN$, $\CN_i$, and $\CC_l,\,l\in\clqiG$, the following relationship holds (\cite{Sakurama2021-oy}):
\begin{equation}\label{eq:neighbor_clique}
    \CN_i=\bigcup_{l\in \clqiG} \CC_l.
\end{equation}


\section{Problem Statement}\label{sec:problem}
Consider a multi-agent system with $n$ agents.
Let $\CN=\{1,\ldots,n\}$ be the set of agent indices.
The communication network is expressed by a time-invariant undirected graph $G=(\CN,\CE)$ with an edge set $\CE$, representing communication paths.
For the graph $G$ and $l\in\clqG=\{1,\ldots,q\}$, $\CC_l$ denotes the $l$th maximal clique of $G$.

Our aim is to design a distributed algorithm that can solve the following distributed optimization problem only using local data and peer-to-peer communication:
\begin{subequations}\label{problem}
\begin{align}
\label{eq:objective_functuon}
    \displaystyle \min_{x\in \mathbb{R}^{nd}} \quad& f(x)= 
    \sum_{i=1}^n f_i(x_i) \\
    \label{eq:constraints}
    \mathrm{s.t.} \quad &x\in \CD = \bigcap_{l\in\clqG} \{x\in\BRnd: x_{\CC_l}\in\CD_l\}
\end{align}
\end{subequations}
with a convex objective function $f$ in \eqref{eq:objective_functuon} and a non-empty closed convex constraint set $\CD$ in \eqref{eq:constraints},
where $x=[x_1^\top,\ldots,x_n^\top]^\top\in\BRnd$.
Here, each $\CD_l\subset\BRCd$ is non-empty and closed convex.

\begin{figure}[t]
    \centering
    \includegraphics[width=0.5\columnwidth]{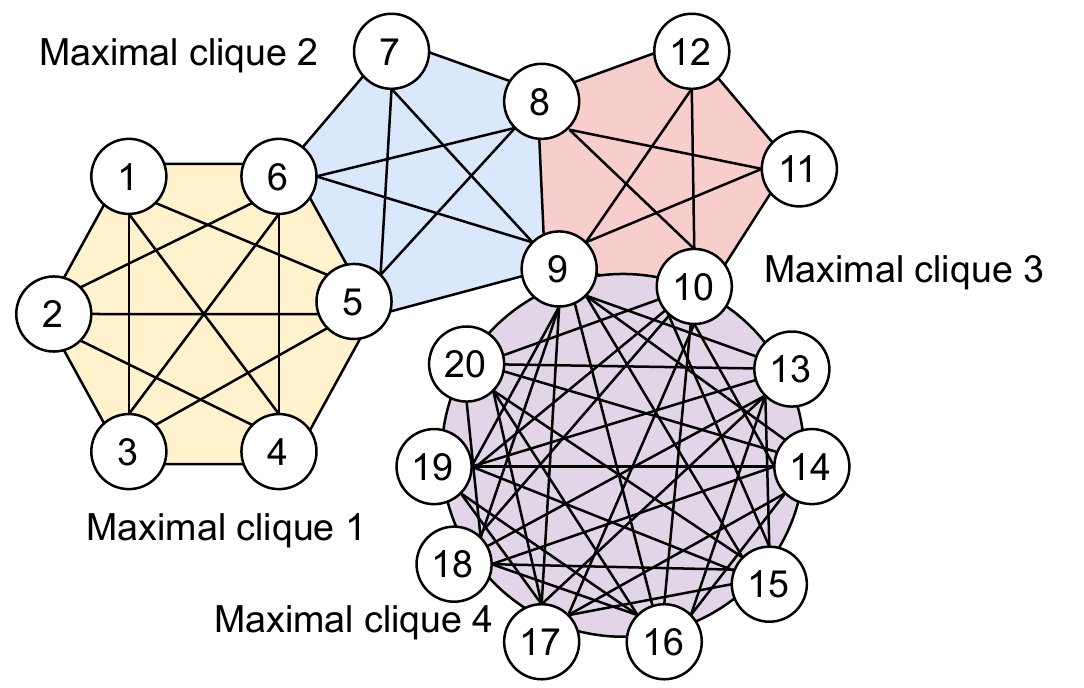} 
    \caption{Example of a network and its maximal cliques.}
    \label{fig:sim_clique}
\end{figure}

The set $\CD$ in \eqref{eq:constraints} can describe various coupled-constraints in accordance with $G$.
The following example represents multiple communities forming a communication network with some constraints.
\begin{exmp}
    Consider graph $G$ in Fig.\ \ref{fig:sim_clique}.
    This network $G$ consists of four communities expressed by maximal cliques, and each of them has some overlapping nodes.
    Consider convex constraints imposed on each communities such as $A_l x_{\CC_l} = b_l$ and/or $\|x_{\CC_l}\| \leq r_l$ for $l\in\clqG=\{1,2,3,4\}$
    with some $A_l\in\mathbb{R}^{m \times |\CC_l|d}$, $b_l\in\mathbb{R}^m$, and $r_l>0$.
    Then, the intersection of these constraints can be written as \eqref{eq:constraints}.
\end{exmp}
Note that conventional agent-wise and pair-wise constraints are included in the class of $\CD$ in \eqref{eq:constraints} as follows.
    \begin{exmp}
    Consider a connected graph $G=(\CN,\CE)$ and non-empty closed convex sets $\CX_{ij}\subset\mathbb{R}^{2d},\,(i,j)\in\CE$. Then, $\CD = \bigcap_{(i,j)\in\CE} \{x\in\BRnd: [x_i^\top,x_j^\top]^\top\in \CX_{ij} \}$
    can be rewritten as \eqref{eq:constraints} because each edge belongs to some maximal clique of $G$.
    Also, the conventional constraints, studied in \cite{Nedic2010-xf}, etc, can be expressed by \eqref{eq:constraints} for $\CD =\bigcap_{(i,j)\in\CE}\{x\in\BRnd:x_i=x_j\}\cap\bigcap_{i=1}^n \{x\in\BRnd: x_i\in \CX_i\}$.
\end{exmp}

Finally, we impose the following assumption on $f$.
\begin{assum}\label{assumption:smoooth}
The function $f$ is $L$-smooth, i.e., $\nabla f$ satisfies $\|\nabla f(x)-\nabla f(y)\| \leq  L \|x-y\|$ with $L>0$ for any $x,y\in\BRnd$.
\end{assum}

\section{Clique-based projected gradient descent}\label{sec:algorithm}
\subsection{Algorithm description}\label{subsec:algorithm}
We develop the clique-based projection and propose a novel distributed algorithm to solve the problem \eqref{problem} using the projection.
Let $x_i(k)\in\mathbb{R}^d$ be the estimate $x_i(k)\in\mathbb{R}^d$ of an solution at each iteration $k$, updated by agent $i$.

The clique-based projection is defined as follows.
\begin{defn}\label{def:clique-based_projection}
    For a non-empty closed convex set $\CD\subset\BRnd$ in \eqref{eq:constraints}, a graph $G$, and its maximal cliques $\CC_l,\,l\in\clqG$, the \textit{clique-based projection} $T:\BRnd\to\BRnd$ of $x\in\BRnd$ onto $\CD$ is defined as
    \begin{equation}
    \label{eq:map_T}
    T(x) = [T_1(x_{\CN_1})^\top,\ldots,T_n(x_{\CN_n})^\top]^\top
\end{equation}
with \begin{equation}
    \label{eq:map_T_i}
    T_i(x_{\CN_i}) = \frac{1}{|\mathrm{clq}_i(G)|} \sum_{l\in\mathrm{clq}_i(G)} 
    \left[ P_{\CD_l} (x_{\CC_l}) \right]_{m_{l,i}}
\end{equation}
for each $i\in \CN$, where $P_{\CD_l}:\BRCd\to\CD_l\subset\BRCd$ is the convex projection of $x_{\CC_l}$ onto $\CD_l$ with respect to the norm $\|\cdot\|_{\mathrm{diag}(\gamma_{\CC_l})}$, i.e.,
\begin{equation}
    \label{eq:map_P_DC}
    P_{\CD_l} (x_{\CC_l})= \argmin_{z\in \CD_l} \|z-x_{\CC_l}\|_{\mathrm{diag}(\gamma_{\CC_l})}.
\end{equation}
 Here, $\gamma=[1/|\mathrm{clq}_1(G)|,\ldots,1/|\mathrm{clq}_n(G)|]^\top \otimes \mathbf{1}_d \in\mathbb{R}^{nd}$, and $m_{l,i}\in \{1,\ldots,|\CC_l|\}$ denotes an order of $i$ in $\CC_l$ for $i\in\CN$ and a clique $\CC_l,\,l\in\clqiG$, i.e., $\CC_l = \{\ldots,\underset{\scriptsize{m_{l,i}}}{i},\ldots\}$.
\end{defn}


Using the clique-based projection, we present a novel distributed algorithm for the problem \eqref{problem}, the \textit{clique-based projected gradient descent} (CPGD), in Algorithm \ref{alg:proposed}.
Here, $x^{[s]}_{\CN_i}(k)$ is the aggregated vector of $x_j^{[s]}(k)$ according to $j\in\CN_i$.
The step size $\{\lambda_k\}$ is a sequence of positive numbers and $p$ is a positive integer.
To guarantee the convergence of the CPGD, we must impose some conditions on $\{\lambda_k\}$.
Section \ref{sec:convergence} shows convergence properties of Algorithm \ref{alg:proposed}.

Note that from \eqref{eq:neighbor_clique}, each agent $i$ can implement the CPGD in Algorithm \ref{alg:proposed} only with local communication with neighbors $\CN_i$.

\begin{figure}[!t]
\begin{algorithm}[H]
\caption{
Clique-based projected gradient descent (CPGD) for agent $i$
}
\label{alg:proposed}
\begin{algorithmic}[1]
\Require $x_i(0)\in\mathbb{R}^d$,\: $\{\lambda_k\}_{k\geq 1}\subset\mathbb{R}$,\:$p\in\mathbb{N}$
\For{$k = 0,1,\ldots$}
\State $x_i^{[1]}(k) \leftarrow x_i(k)- \lambda_{k+1} \nabla_i f_i(x_i(k))$.
\For{$s = 1,2,\ldots p$}
\State Gather $x_j^{[s]}(k)$ from each $j\in \CN_i\setminus \{i\}$.
\State $x_i^{[s+1]}(k) \leftarrow T_i(x_{\CN_i}^{[s]}(k))$ with $T_i$ in \eqref{eq:map_T_i}.
\EndFor
\State $x_i(k+1) \leftarrow x_{i}^{[p+1]}(k).$
\EndFor
\end{algorithmic}
\end{algorithm}
\end{figure}

By aggregating Algorithm \ref{alg:proposed} for all $i\in\CN$, we obtain
\begin{subequations}
\label{eq:proposed_method_full}
\begin{align}
    x^{[1]}(k)&=x(k) - \lambda_{k+1} \nabla f(x(k)) \\
    x(k+1)    &= T^p(x^{[1]}(k)),
\end{align}
\end{subequations}
where $T^p=\overbrace{T\circ T\circ \cdots \circ T}^{p}$.

\subsection{Discussion on the clique-based projection}\label{subsec:CP_analysis}
Here, we view the clique-based projection in Definition \ref{def:clique-based_projection} from the perspective of the convex projection and the proximal operator.

First, we prove the following proposition, which indicates that the clique-based projection $T$ can be obtained by taking the gradient of a function $V$.
\begin{prop}\label{prop:V_T}
    Let $V:\BRnd\to\mathbb{R}$ be
\begin{equation}\label{eq:V}
    V(x) = \frac{1}{2} \sum_{l\in\clqG} \|x_{\CC_l}-P_{\CD_l}(x_{\CC_l})\|^2_{\mathrm{diag}(\gamma_{\CC_l})}
\end{equation}
with $P_{\CD_l}$ in \eqref{eq:map_P_DC} and $\gamma=[1/|\mathrm{clq}_1(G)|,\ldots,1/|\mathrm{clq}_n(G)|]^\top \otimes \mathbf{1}_d \in\mathbb{R}^{nd}$. Then, the following holds for any $x\in\BRnd$:
\begin{equation}\label{eq:V_T}
    T(x) = x - \nabla V(x).
\end{equation}
\end{prop}
 \begin{pf}
     Since each $\CD_l$ is closed and convex, $1/2\,\|x_{\CC_l}-P_{\CD_l}(x_{\CC_l}) \|_{\mathrm{diag}(\gamma_{\CC_l})}^2$ is differentiable and thus $V(x)$ in \eqref{eq:V} is also differentiable.
     Then, for all $i\in\CN$, we have
        $\nabla_i V(x) 
          = \sum_{l\in \clqiG} \frac{1}{|\clqiG|} (x_i - [P_{\CD_l}(x_{\CC_l})]_{m_{l.i}})
         = x_i - \sum_{l\in \clqiG} \frac{1}{|\clqiG|} [P_{\CD_l}(x_{\CC_l})]_{m_{l.i}}
         = x_i -T_i(x)$
     from \eqref{eq:neighbor_clique} and \eqref{eq:map_T_i}.
     Hence, we obtain \eqref{eq:V_T}. $\square$
 \end{pf}

    Through Proposition \ref{prop:V_T}, we can interpret the proposed method in \eqref{eq:map_T} as a variant of the proximal gradient method (\cite{Beck2009-bq}).
    From \eqref{eq:V_T}, the clique-based projection $T$ satisfies 
    \begin{equation*}
        T(x) = \argmin_{y\in\BRnd} \: \frac{1}{2}\|x-y\|^2 + V(x) + \nabla V(x)^\top(y-x)
        .
    \end{equation*}
    Hence, the clique-based projection $T$ can be regarded as the proximal operator $\mathrm{prox}_{\psi}(x) = \argmin_{y\in\BRnd} 1/2 \|x-y\|^2 +\psi(y)$ for $\psi(y) = V(x) + \nabla V(x)^\top(y-x) $, which is
    the first order approximation of $V(y)$ at $x$.
    
    Moreover, if $G$ is complete, the CPGD in Algorithm \ref{alg:proposed} equals to the well-known \textit{projected gradient descent} (PGD) (see \cite{Calamai1987-ai}):
    \begin{equation}\label{eq:PGD}
        x(k+1) = P_\CD (x(k) - \lambda_{k+1} \nabla f(x(k)))
    \end{equation}
    with $P_\CD(x)=\argmin_{z\in\CD}\|x-z\|$
    because $\clqG=\{1\}$ and $\CC_1 = \CN$ hold for complete $G$.
    Note that the PGD in \eqref{eq:PGD} is centralized due to the operation $P_{\CD}(\cdot)$.
\section{Convergence Analysis}\label{sec:convergence}

For the CPGD in Algorithm \ref{alg:proposed}, we present convergence theorems for both diminishing and fixed step sizes.

Before proceeding to analyze convergence properties, we provide the following lemma, which shows key features of the clique-based projection $T$ in Definition \ref{def:clique-based_projection}.
\begin{lem}\label{lemma:T}
    For the clique-based projection $T$ in Definition \ref{def:clique-based_projection} and the closed convex set $\CD$ in \eqref{eq:constraints}, the following statements hold:
    \begin{itemize}
    \setlength{\itemsep}{0.05cm} 
        \item[a)] The mapping $T$ is nonexpansive, i.e., $\|T(x)-T(y)\|\leq \|x-y\|$ holds for any $x,y\in\BRnd$.
        \item[b)] The fixed points set of $T$ satisfies $\mathrm{Fix}(T)=\CD$.
        \item[c)] For any $x\in\BRnd\setminus\CD$ and any $z\in\CD$, $\|T(x)-z\|<\|x-z\|$ holds.
        \item[d)] For any $x\in\BRnd$, $T^\infty(x)=\lim_{p\to \infty} T^p(x) \in \CD$ holds.
    \end{itemize}
\end{lem}
\begin{pf}
    See Appendix \ref{sec:proof}. $\square$
\end{pf}
By Lemma \ref{lemma:T}a-b, the clique-based projection $T$ preserves important features of the convex projection for $\CD$, that is, the nonexpansiveness and fixed point set.
Moreover, from Lemma \ref{lemma:T}c-d, by repeatedly operating $T$, 
$x$ gradually converges to $\CD$.

With this in mind, we show the first main result as follows.
\begin{thm}\label{theorem:non_acc}
Assume that a convex objective function $f$ in \eqref{eq:objective_functuon} satisfies Assumption \ref{assumption:smoooth}, and that $\CD\subset\BRnd$ in \eqref{eq:constraints} is a non-empty closed convex set.
Consider the sequence $\{x(k)\}$ generated by CPGD in Algorithm \ref{alg:proposed} (or \eqref{eq:proposed_method_full}).
\begin{itemize}
    \item[a)] Let a sequence $\{\lambda_k\}$ of positive integers satisfy $\lim_{k\to \infty} \lambda_k=0$, $\sum_{k=1}^\infty \lambda_k =\infty$, and $\sum_{k=1}^\infty \lambda_k^2 < \infty$.\footnote{For example, $\lambda_k=1/k$ satisfies the conditions.}
    Assume that $\CD$ in \eqref{eq:constraints} is bounded.
    Then, for any initial point $x(0)=x_0\in\BRnd$ and any $p\in\mathbb{N}$, $\{x(k)\}$ converges to an optimal solution $x_*\in\argmin_{x\in\CD}f(x)$.
     \item[b)] Let a sequence $\{\lambda_k\}$ of positive integers satisfy $\lim_{k\to \infty} \lambda_k=0$, $\sum_{k=1}^\infty \lambda_k =\infty$, and $\sum_{k=1}^\infty |\lambda_k - \lambda_{k+1}| < \infty$.
    \footnote{For example, $\lambda_k=1/k$ and $\lambda_k=1/\sqrt{k}$ satisfy the conditions.}
    Additionally assume that the gradient $\nabla f:\BRnd\to\BRnd$ of $f$ is strongly monotone, i.e., there exists some $ \mu>0$ such that
    $(\nabla f(x)-\nabla f(y))^\top (x-y) \geq  \mu \|x-y\|^2 $ is satisfied for any $x,y\in\BRnd$.
 Then $\{x(k)\}$ converges to the unique optimal solution $x_*=\argmin_{z\in\CD} f(z)$ for any initial point $x(0)=x_0\in\BRnd$ and any $p\in\mathbb{N}$.
    \item[c)] Let $\lambda_k=t\in(0,1/L]$ for any $k\in\mathbb{N}$.
    Let $J:\BRnd\to\mathbb{R}$ be 
    \begin{equation}\label{eq:J}
    J(x)=f(x)+ V(x)/t    
    \end{equation}
    with $V$ in \eqref{eq:V}.
    Then, for any initial point $x(0)=x_0\in\BRnd$ and $p=1$,
    \begin{equation}\label{eq:fixed_bound}
         J(x(k))-J(x_*)\leq \frac{\|x_0-x_*\|^2}{2t k} 
    \end{equation}
    holds with $x_*\in \argmin_{x\in\CD} f(x)$.
\end{itemize}
\end{thm}
\begin{pf}
    a) From Lemma \ref{lemma:T}a-b, the CPGD in \eqref{eq:proposed_method_full} can be regarded as the hybrid steepest descent in \cite{Yamada2001-xv,Yamada2002-pc} for any $p\in\mathbb{N}$.
    Hence, Theorem \ref{theorem:non_acc}a follows from Theorem 2.18, Remark 2.17 in \cite{Yamada2002-pc}, and Lemma \ref{lemma:T}c.
    b) The statement follows from Theorem 2.15 in \cite{Yamada2002-pc} and Lemma \ref{lemma:T}a-b.
    c) See Appendix \ref{subsec:supporting_lemma} and \ref{subsec:proof_fixed}. $\square$
\end{pf}

Theorem \ref{theorem:non_acc}a-b imply that the CPGD with a diminishing step size provides an optimal solution regardless of the number $p$ of operations of $T$ in each iteration.
Thus, with a sufficiently large $p$, the CPGD can generate points arbitrarily close to $\CD$ by Lemma \ref{lemma:T}c-d and thus can work just like the conventional PGD in \eqref{eq:PGD}, which is centralized.
Note that the assumption of the boundedness of $\CD$ in Theorem \ref{theorem:non_acc}a is not restrictive in practice, and the additional assumption of the strong monotonicity of $\nabla f$ in Theorem \ref{theorem:non_acc}b is satisfied for strongly convex $f$.

Theorem \ref{theorem:non_acc}c guarantees the non-ergodic convergence rate of $O(1/k)$ for smooth $f$ when $\lambda_k$ is fixed and $p=1$.
Practically, the CPGD is expected to perform well even if $p>1$, which is shown in the numerical results in Section \ref{sec:numerical_experiments}.


\section{Nesterov's Accelerated CPGD}\label{sec:acc}
In this section, to enhance the convergence rate, we modify the proposed CPGD with Nesterov's acceleration scheme (\cite{Nesterov1983-ir,Beck2009-bq}).
The modified method is called the accelerated clique-based projected gradient descent (ACPGD), given as follows:
\begin{subequations}
\label{eq:nesterov}
\begin{align}
     x(k+1)&=T^p(\hat{x}(k) - \lambda_{k+1} \nabla f(\hat{x}(k)) )\\
  \hat{x}(k+1)&= x(k+1) + \frac{\sigma_k-1}{\sigma_{k+1}} (x(k+1)-x(k)),
\end{align}
\end{subequations}
where $\hat{x}(0)=x(0)$ and $\sigma_{k+1}=(1+\sqrt{1+4\sigma_k^2})/2,\,\sigma_0=1$. 
Agent $i\in\CN$ can run the ACPGD in a distributed fashion by updating $x_i(k+1)$ as $\hat{x}_i(k+1)= x_i(k+1) + \frac{\sigma_k-1}{\sigma_{k+1}} (x_i(k+1)-x_i(k))$ with $x_i(k)$ in addition to Algorithm \ref{alg:proposed}.

The ACPGD achieves the convergence rate of $O(1/k^2)$ as follows.
\begin{thm}\label{theorem:acceleration}
    Assume that a convex function $f$ in \eqref{eq:objective_functuon} satisfies Assumption \ref{assumption:smoooth}, and that $\CD\subset\BRnd$ in \eqref{eq:constraints} is a non-empty closed convex set.
    Let $p=1$ and $\lambda_k=t \in (0,1/ L]$ for all $k\in\mathbb{N}$.
    Consider the sequence $\{x(k)\}$ generated by the ACPGD in \eqref{eq:nesterov}.
    Then, for any initial state $x(0)=\hat{x}(0)=x_0\in\BRnd$, the following inequality holds:
    \begin{equation}\label{eq:nesterov_bound}
        J(x(k)) - J(x_*) \leq \frac{2\|x_0-x_*\|^2}{tk^2},
    \end{equation}
     where $x_* \in \argmin_{x\in\CD} f(x)$ and $J(x)$ is given as \eqref{eq:J}.
\end{thm}
\begin{pf}
See Appendix \ref{subsec:supporting_lemma} and \ref{subsec:proof_acc}. $\square$
\end{pf}

Although the case of $p=1$ is only proved in Theorem \ref{theorem:acceleration} like Theorem \ref{theorem:non_acc}c, 
the ACPGD for $p>1$ can also perform well as shown in numerical results.

\section{Numerical Experiment}\label{sec:numerical_experiments}

We demonstrate the effectiveness of the proposed method through numerical experiments.

Consider a multi-agent system with $n=20$ agents.
The communication network $G$ is given as Fig.\ \ref{fig:sim_clique}.
Then we have $\clqG=\{1,2,3,4\}$ and $\CC_1=\{1,2,\ldots,6\},\,\CC_2=\{5,6,\ldots,9\},\,\CC_3=\{8,9,\ldots,12\},\,\CC_4=\{9,10,13,14,\ldots,20\}$.
We consider the following allocation problem:
\begin{subequations}
\label{eq:num_problem}
\begin{align}
    \min_{x\in\mathbb{R}^{20}} &\quad \frac{1}{2}\sum_{i=1}^{20} (x_i-a_i)^2 \label{eq:num_objective}\\
    \mathrm{s.t.}&\quad \sum_{j\in\CC_l} x_j = N_l,\; \forall l\in\clqG=\{1,\ldots,4\},\label{eq:num_constraints1}
\end{align}
\end{subequations}
where $a_1,\ldots,a_{20} \in\mathbb{R}$ are randomly generated by the uniform distribution for the interval $[0,10]$ and $N_l,\,l\in\clqG$ are given as $[N_1,\,\ldots,N_4] = [7,3,5,10]$.
Letting $\CD_l=\{y=[y_1,\ldots,y_{|\CC_l|}]^\top\in\mathbb{R}^{|\CC_l|}: \sum_{j=1}^{|\CC_l|} y_j = N_l \}$ for $l\in\clqG$, we can apply Algorithm \ref{alg:proposed} to the problem \eqref{eq:num_problem}.

We conduct simulations of the CPGD in Algorithm \ref{alg:proposed} with  $\lambda_k=1/k$ and $\lambda_k=0.001$ for $p=1,\,10,\,50$.
Additionally, we conduct simulations of the ACPGD in \eqref{eq:nesterov} with $\lambda_k=0.001$ for $p=1,10,50$.
For comparison, we run the conventional PGD in \eqref{eq:PGD}, which is centralized.

\begin{figure*}[t]
    \centering
    \hspace*{-0.25cm}\includegraphics[width=2\columnwidth]{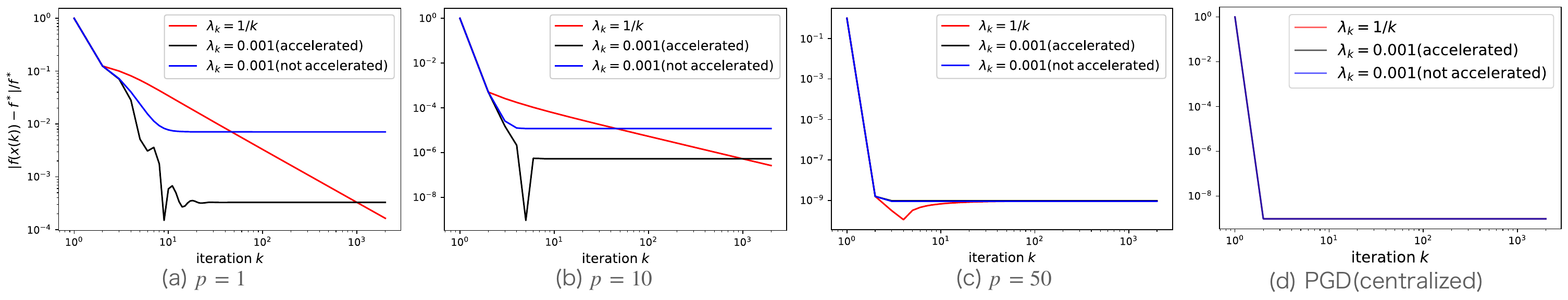} 
    \caption{Plots of relative residuals of $f$ versus the number $k$ of iterations under the CPGD with
    $\lambda_k=1/k$ (red line) and $\lambda_k=0.001$ (blue line), and the ACPGD with $\lambda_k=0.001$ (black line) for (a) $p=1$, (b) $p=10$, and (c) $p=50$.
    For comparison purposes,
    the plot (d) shows the results of the PGD in \eqref{eq:PGD} with $\lambda_k=1/k,\,0.001$ and the accelerated PGD with $\lambda_k=0.001$.
    }
    \label{fig:residual}
\end{figure*}

Figs.\ \ref{fig:residual}a-\ref{fig:residual}c plot the evolution of the relative optimality gap $|f(x(k))-f^*|/f^*$ between the value of the objective function $f(x(k))$ and its optimal value $f^*$ under $p=1,\,10,\,50$, respectively.
Besides, Fig.\ \ref{fig:residual}d shows the result of the conventional PGD.
From Figs.\ \ref{fig:residual}a-\ref{fig:residual}c, the CPGD with both types of step sizes and the ACPGD succeed in the swift convergence of $x(k)$ to a point close to the optimal solution.
For a larger $p$, a more accurate solution is obtained.
In addition, in the early stage of the simulation, the ACPGD achieves the best performance for all $p$, especially the lower $p$.
After about $1000$ iterations, the diminishing step size $\lambda_k=1/k$ gives the closest solution to the optimal one.
Furthermore, from Fig.\ \ref{fig:residual}c and Fig.\ \ref{fig:residual}d, the CPGD and the ACPGD run as fast as the PGD when $p=50$.
These results illustrate the effectiveness of the proposed method.

\section{Conclusion}\label{sec:conclusion}
This paper addressed a distributed convex optimization problem with clique-wise couplings in constraints.
First, we developed the clique-based projection operator and proposed a new distributed algorithm with the operator,
the clique-based projected gradient descent (CPGD).
Next, we proved its convergence properties for diminishing and fixed step sizes.
Moreover, we presented the accelerated version of the CPGD, which achieved the convergence rate of $O(1/k^2)$.
Finally, numerical experiments illustrated the effectiveness of the proposed method.
Our future directions are to consider more general coupled constraints and to apply the CPGD to some applications, such as traffic networks.

\appendix
\section{Proof of Lemma \ref{lemma:T}}\label{sec:proof}
\subsection{Preliminaries}\label{subsec:proof_pre}
As a preliminary, we provide several useful inequalities and some propositions.

First, for a convex function $h:\mathbb{R}^m\to\mathbb{R}$, any positive integer $k\in\mathbb{N}$, any $x_1,\ldots,x_k\in\mathbb{R}^m$, and any $\alpha_j\geq 0\;(j=1,\ldots,k)$ satisfying $\sum_{j=1}^k \alpha_j=1$, the following inequality holds, which is called \textit{Jensen's inequality}:
\begin{equation}\label{Jensen}
    h(\textstyle\sum_{j=1}^k \alpha_j x_j ) \leq \textstyle\sum_{j=1}^k \alpha_j h(x_j).
\end{equation}
Note that $h( \sum_{j=1}^k \alpha_j x_j ) = \sum_{j=1}^k \alpha_j h(x_j)$ holds if and only if $x_1=\cdots=x_k$ holds or $h$ is affine (see \cite{Peressini1988-ak}).

Second, we give some properties of the convex projection without proof.
For details, see textbooks of function analysis, e.g., \cite{Kreyszig1991-rz}.
For any $x,y\in\mathbb{R}^m$, a norm $\|\cdot\|_Q$, a non-empty closed convex set $\mathcal{M} \subset\mathbb{R}^m$, and any $z\in\mathcal{M}$, the convex projection $P_{\mathcal{M}}:\mathbb{R}^m\to\mathcal{M}$ with respect to $\|\cdot\|_Q$, i.e., $P_\mathcal{M}(x)=\argmin_{z\in\mathcal{M}}\|x-z\|_Q$, satisfies the following inequalities:
\begin{align}
\label{nonexpansive_P}
    \|P_{\mathcal{M}}(x)-P_{\mathcal{M}}(y)\|_Q &\leq \|x-y\|_Q \\
\label{quasi-nonexpansive_P}
    \|P_{\mathcal{M}}(x)-z\|_Q &\leq \|x-z\|_Q,
\end{align}
implying the nonexpansiveness and quasi-nonexpansiveness of $P_\mathcal{M}(\cdot)$, respectively. 
 Moreover, for any $x,y\in\mathrm{R}^m$, the following inequalities hold:
 \begin{align}\label{eq:P_ip}
    &\!\!\!\!
    \|P_{\mathcal{M}}(x)-P_{\mathcal{M}}(y) \|_Q^2 \leq (x-y)^\top Q (P_{\mathcal{M}}(x)-P_{\mathcal{M}}(y))  \!\!\\
     \label{eq:P_equiv}
     &(x-P_{\mathcal{M}}(x))^\top Q (z-P_{\mathcal{M}}(x)) \leq 0\quad (\forall z\in\mathcal{M}).
 \end{align}

Next, we present important properties of the function $V(x)$ in \eqref{eq:V} for $\CD$ in \eqref{eq:constraints} as follows. Note that the function $V$ in \eqref{eq:V} is convex because of the convexity of each $\CD_l$.
\begin{prop}\label{prop:V}
For $V(x)$ in \eqref{eq:V} and a non-empty closed convex set $\CD$ in \eqref{eq:constraints}, $V(x)=0 \Leftrightarrow x\in \CD$ holds.
\end{prop}
\begin{pf}
    If $V(x)=0$ for $x\in\BRnd$, we obtain $x_{\CC_l}=P_{\CD_l}(x_{\CC_l}) \in \CD_l$ for all $l\in\clqG$, which yields $x\in\CD$ because of \eqref{eq:constraints}.
    Conversely, if $x\in\CD$, then we have $x_{\CC_l}\in \CD_l$ for all $l\in\clqG$.
    Thus, $V(x)=0$ holds. $\square$
\end{pf}
\begin{prop}\label{prop:V_smooth}
    The function $V(x)$ in \eqref{eq:V} is a $1$-smooth function, i.e., its gradient $\nabla V(x)$ is $1$-Lipschitzian.
\end{prop}
\begin{pf}
    From Definition \ref{def:clique-based_projection}, the inequality \eqref{eq:P_ip}, and Proposition \ref{prop:V_T}, we obtain the following for any $x,y\in\BRnd$:
    \begin{align*}
        &\|\nabla V(x)- \nabla V(y)\|^2 = \|(x-y) - (T(x)-T(y))\|^2 \\
        =& \|x-y\|^2 + \|T(x)-T(y)\|^2 -2 (x-y)^\top(T(x)-T(y)) \\
        =&\|x-y\|^2 + \|T(x)-T(y)\|^2 \\
         & -2 \textstyle\sum_{l\in\clqG} (x_{\CC_l}-y_{\CC_l})^\top \mathrm{diag}(\gamma_{\CC_l})(P_{\CD_l}(x_{\CC_l})-P_{\CD_l}(y_{\CC_l})) \\
        \leq & \|x-y\|^2+\|T(x)-T(y)\|^2\\
         &\quad -2 \textstyle\sum_{l\in\clqG}  \|P_{\CD_l}(x_{\CC_l})-P_{\CD_l}(y_{\CC_l})\|^2_{\mathrm{diag}(\gamma_{\CC_l})} \\
         \leq & \|x-y\|^2-\|T(x)-T(y)\|^2\leq\|x-y\|^2.
     \end{align*}
     The last line follows from \eqref{eq:T_nonexpansive_proof} in the proof of Lemma \ref{lemma:T}a.
     It completes the proof. $\square$
\end{pf}
\subsection{Proof}\label{subsec:proof_theorem_T}
Here, we show the proof of Lemma \ref{lemma:T}.

a) 
From \eqref{Jensen} and
\eqref{nonexpansive_P}, the following inequality holds:
\begin{align}\label{eq:T_nonexpansive_proof}
    &\|x-y\|^2
    = \textstyle\sum_{l \in \clqG} \|x_{\CC_l}-y_{\CC_l}\|_{\mathrm{diag}(\gamma_{\CC_l})}^2 \nonumber \\
     \geq& \textstyle\sum_{l\in \clqG} \|P_{\CD_l} (x_{\CC_l})-P_{\CD_l} (y_{\CC_l})\|^2_{\mathrm{diag}(\gamma_{\CC_l})} \nonumber \\
     =& \textstyle\sum_{i=1}^n\textstyle\sum_{l\in \clqiG} \frac{1}{|\clqiG|} \|\left[ P_{\CD_l} (x_{\CC_l}) \right]_{m_{l,i}} \nonumber \\
      &\quad\quad \quad\quad \quad\quad \quad\quad\quad\quad \quad\quad\quad
      -\left[ P_{\CD_l} (y_{\CC_l}) \right]_{m_{l,i}}\|^2 \nonumber\\
     \geq& \textstyle\sum_{i=1}^n \|\sum_{l\in \clqiG} \frac{1}{|\clqiG|}(\left[ P_{\CD_l} (x_{\CC_l}) \right]_{m_{l,i}} 
     \nonumber \\
     &\quad\quad \quad\quad \quad\quad \quad\quad\quad\quad \quad\quad\quad
     -\left[ P_{\CD_l} (y_{\CC_l}) \right]_{m_{l,i}}) \|^2
     \nonumber \\
     =& \textstyle\sum_{i=1}^n \|T_i(x_{\CN_i}) - T_i(y_{\CN_i})\|^2 
     = \|T(x) - T(y)\|^2.
\end{align}
Thus, we obtain $\|T(x)-T(y)\|\leq \|x-y\|$. $\square$

b) $\CD\subset \mathrm{Fix}(T)$ holds because $x_{\CC_l}=P_{\CD_l}(x_{\CC_l})$ holds for any $x\in\CD$ and all $l\in\clqG$.
In the following, we prove the converse inclusion $\mathrm{Fix}(T)\subset\CD$.
Let $y\in\CD$. We show $\hat{y}\in \mathrm{Fix}(T)\setminus\{y\}  \Rightarrow \hat{y}\in \CD$.
From $\hat{y}\in\mathrm{Fix}(T)$, we obtain $\hat{y}_i=T_i(\hat{y}_{\CN_i})$ for all $i\in\CN$.
In addition, from \eqref{Jensen} and \eqref{quasi-nonexpansive_P}, we have
\begin{align*}
\small
    &\|y-\hat{y}\|^2
    \geq \textstyle\sum_{l\in\clqG} \|y_{\CC_l} - P_{\CD_{\CC_l}}(\hat{y}_\CC)\|^2_{\mathrm{diag}(\gamma_{\CC_l})} \nonumber\\
    &=\textstyle\sum_{i=1}^n \sum_{l\in \mathrm{clq}_i(G)} \frac{1}{|\mathrm{clq}_i(G)|} \|y_i-\left[ P_{\CD_l} (\hat{y}_{\CC_l}) \right]_{m_{l,i}}\|^2  \nonumber \\
   &\geq \textstyle\sum_{i=1}^n  \| y_i - \textstyle\sum_{l\in \mathrm{clq}_i(G)} \frac{1}{|\mathrm{clq}_i(G)|} \left[ P_{\CD_l} (\hat{y}_{\CC_l}) \right]_{m_{l,i}}\|^2
    \nonumber \\
   &= \textstyle\sum_{i=1}^n  \| y_i - T_i(\hat{y}_{\CN_i})\|^2 = \|y-\hat{y}\|^2.
\end{align*}
Thus, from the equality condition of \eqref{Jensen}, we obtain $y_i-\left[ P_{\CD_k} (\hat{y}_{\CC_k}) \right]_{m_{k,i}}=y_i-\left[ P_{\CD_l} (\hat{y}_{\CC_l}) \right]_{m_{l,i}}$ for all $\CC_k,\CC_l\,(k,l\in\clqiG)$ for all $i\in\CN$.
Then, we have $\left[ P_{\CD_k} (\hat{y}_{\CC_k}) \right]_{m_{k,i}}=\left[ P_{\CD_l} (\hat{y}_{\CC_l}) \right]_{m_{l,i}}$ for all $\CC_k,\CC_l\,(k,l\in\clqiG)$.
Therefore, because $\hat{y}\in\mathrm{Fix}(T)$, we have $2V(\hat{y})=\sum_{i=1}^n \sum_{l\in \mathrm{clq}_i(G)} \frac{1}{|\mathrm{clq}_i(G)|} \|\hat{y}_i-\left[ P_{\CD_l} (\hat{y}_{\CC_l}) \right]_{m_{l,i}}\|^2= \sum_{i=1}^n \| \hat{y}_i - T_i(\hat{y}_{\CN_i})\|^2=0$.
Then, $\hat{y}\in\CD$ holds from Proposition \ref{prop:V}.
Hence, we obtain $\mathrm{Fix}(T)\subset\CD$.
$\square$

c) For a non-empty closed convex set $\CD$ in \eqref{eq:constraints} and $x\in\BRnd\setminus\CD$, there exists $\hl\in\clqG$ such that $\|x_{\CC_{\hl}}-P_{\CD_{\hl}}(x_{\CC_{\hl}})\|_{\mathrm{diag}(\gamma_{\CC_{\hl}})}>0$.
Hence, for $\hl\in\clqG$, $x\in\BRnd\setminus\CD$, and $z\in\CD$, we have $\|x_{\CC_{\hl}}-z_{\CC_{\hl}}\|^2_{\mathrm{diag}(\gamma_{\CC_{\hl}})}>\|P_{\CD_{\hl}}(x_{\CC_{\hl}})-z_{\CC_{\hl}}\|^2_{\mathrm{diag}(\gamma_{\CC_{\hl}})}$ because $\|x_{\CC_{\hl}}-z_{\CC_{\hl}}\|^2_{\mathrm{diag}(\gamma_{\CC_{\hl}})} 
    = \|x_{\CC_{\hl}}-P_{\CD_{\hl}}(x_{\CC_{\hl}})\|^2_{\mathrm{diag}(\gamma_{\CC_{\hl}})} + \|P_{\CD_{\hl}}(x_{\CC_{\hl}})-z_{\CC_{\hl}}\|^2_{\mathrm{diag}(\gamma_{\CC_{\hl}})} 
     -2 (x_{\CC_{\hl}}-P_{\CD_{\hl}}(x_{\CC_{\hl}}))^\top \mathrm{diag}(\gamma_{\CC_{\hl}})(z_{\CC_{\hl}}-P_{\CD_{\hl}}(x_{\CC_{\hl}})) 
    > \|P_{\CD_{\hl}}(x_{\CC_{\hl}})-z_{\CC_{\hl}}\|^2_{\mathrm{diag}(\gamma_{\CC_{\hl}})}$ holds,
where the last line follows from \eqref{eq:P_equiv}.
Thus, by \eqref{Jensen} and \eqref{nonexpansive_P}, for any $x\in\BRnd\setminus\CD$ and $z\in\CD$,
we obtain
\begin{align*}
    \|x-z\|^2 =& \textstyle\sum_{l\in \clqG} \|x_{\CC_l}-z_{\CC_l}\|^2_{\mathrm{diag}(\gamma_{\CC_l})}  \\
     >& \textstyle\sum_{l\in \clqG} \|P_{\CD_l}(x_{\CC_l})-z_{\CC_l}\|^2_{\mathrm{diag}(\gamma_{\CC_l})}  \\
    \geq& \textstyle\sum_{i=1}^n \|\textstyle\sum_{l\in\clqiG} \frac{1}{|\clqiG|}[P_{\CD_l} (x_{\CC_l})]_{m_{l,i}}-z_i \|^2 \\
    =& \textstyle\sum_{i=1}^n \|T_i(x_{\CN_i}) -z_i\|^2 = \|T(x)-z\|^2.
\end{align*}
Therefore, $\|T(x)-z\|<\|x-z\|$ holds for any $x\in\BRnd\setminus\CD$ and any $z\in\CD$.
$\square$

d) For $x\in\BRnd$, we define $\{a_k\}$ as $a_{k+1} = T(a_k)$ with $a_0=x$.
Then, we obtain $\lim_{k\to\infty} a_{k+1} = \lim_{k\to\infty} T(a_k)$.
Thus, from the continuity of $T$ shown in Lemma \ref{lemma:T}a, we have $T^\infty(x)= \lim_{k\to\infty} a_{k+1} =T( \lim_{k\to\infty} a_k) = T(T^\infty(x))$.
Hence, Lemma \ref{lemma:T}b yields $T^\infty (x)\in\mathrm{Fix}(T) = \CD$. $\square$

\section{Proof of Theorem \ref{theorem:non_acc}c and \ref{theorem:acceleration}}\label{sec:proof_theorem_convergence}
Here, we show the proofs of Theorem \ref{theorem:non_acc}c and \ref{theorem:acceleration}.
These proofs are based on the convergence theorems for ISTA and FISTA (Theorem 3.1 and 4.4 in \cite{Beck2009-bq}), respectively.

\subsection{Supporting lemmas}\label{subsec:supporting_lemma}
Before proceeding to prove the theorems, we show some inequalities corresponding to those obtained from Lemma 2.3 in \cite{Beck2009-bq}, which is a key to prove the convergence theorems.
Note that a differentiable function $h:\mathbb{R}^m\to\mathbb{R}$ is convex if and only if
\begin{equation}\label{eq:convex_inquality}
    h(y)\geq h(x) + \nabla h(x)^\top(y-x)
\end{equation}
holds for any $x,y\in\BRnd$.
Besides, if $h$ is $\beta$-smooth and convex, then
\begin{align}
\label{eq:smooth1}
    h(y) &\leq h(x) + \nabla h(x)^\top(y-x) + \frac{\beta}{2} \|y-x\|^2 \\
\label{eq:smooth2}
    h(y) &\geq h(x) + \nabla h(x)^\top(y-x) + \frac{1}{2\beta} \|\nabla h(x)-\nabla h(y)\|
\end{align}
holds for any $x,y\in\BRnd$. For details, see textbooks on nonlinear programming, e.g.,
\cite{Bertsekas1999-na}.

In preparation for showing lemmas, let $t\in (0,1/L]$ and
\begin{equation*}
V_t(x) = V(x)/t
\end{equation*}
with $V(x)$ in \eqref{eq:V}.
Additionally, for $s\in\BRnd$, we define $F_y :\BRnd\to \mathbb{R}$ with some $y\in\BRnd$ as
\begin{equation}\label{def:F_y}
 F_{y} (s) = f(s) + V_t(y) + \nabla V_t(y)^\top (s-y).
\end{equation}
For $F_y(s)$ in \eqref{def:F_y}, the following inequalities hold.
\begin{prop}
Assume that the convex function $f$ satisfies Assumption \ref{assumption:smoooth}.
    Let $y = x-t\nabla f(x)$ for $x\in \BRnd$.
    Then, 
    \begin{equation}\label{eq:lem_acc_1}
        F_y(T(y)) \leq F_y (z) + \frac{1}{t} (x-T(y))^\top (x-z) - \frac{1}{2t} \|x-T(y)\|^2
    \end{equation}
    holds for any $z\in\BRnd$.
\end{prop}
\begin{pf}
    Let $G_y(s) = f(s) + \nabla V_t(y)^\top (s-y) $ and $z\in\BRnd$. Then, by using $L$-smoothness of $f$ in Assumption \ref{assumption:smoooth}, $\nabla f(x) = (x-y)/t$, and $\nabla V_t (y) = (y-T(y))/t$ (see Proposition \ref{prop:V_T}), 
    \fontsize{9.8pt}{2.5pt}\selectfont
    \begin{align*}
        &G_y(T(y)) 
        = f(T(y)) + \nabla V_t(y)^\top (T(y)-y) \\
        \leq& f(x) - \nabla f(x)^\top (x-T(y)) + \frac{1}{2t} \|x-T(y)\|\\
        &+ \nabla V_t(y)^\top (T(y)-y) \\
        \leq& f(z) + \nabla f(x)^\top (x-z) - \nabla f(x)^\top(x-T(y)) \\
        +& \frac{1}{2t} \|x-T(y)\|^2 + \nabla V_t(y)^\top (z-y) + \nabla V_t(y)^\top (T(y)-z) \\
        =& G_y(z) + \frac{1}{t}(x-T(y))(T(y)-z) + \frac{1}{2t}\|x-T(y)\|^2\\
        =&  G_y(z) + \frac{1}{t}(x-T(y))^\top (x-z) -\frac{1}{2t} \|x-T(y)\| 
        \end{align*}
        \normalsize
        is obtained from \eqref{eq:convex_inquality} and \eqref{eq:smooth1}.
        Thus, adding $V_t(y)$ to the both sides, we obtain \eqref{eq:lem_acc_1}. $\square$
\end{pf}

\begin{prop}
    Let $x(k+1) = T(y(k))$ with some $\{y(k)\} \subset \BRnd$. Then, it holds that
    \begin{align}\label{eq:lem_acc_2}
        &F_{y(k)} (x(k)) + \frac{t}{2} \|\nabla V_t(y(k))\|^2 \notag \\
        \leq& F_{y(k-1)} (x(k)) + \frac{t}{2} \|\nabla V_t(y(k-1))\|^2.
    \end{align}
\end{prop}
\begin{pf}
By $1/t$-smoothness of $V_t(x)$ (see Proposition \ref{prop:V_smooth}) and Proposition \ref{prop:V_T},
\begin{align*}
    &F_{y(k-1)} (x(k)) = f(x(k)) + V_t(y(k-1)) \\
    &+ \nabla V_t(y(k-1))^\top (x(k)-y(k-1)) \\
     =& f(x(k)) + V_t(y(k-1)) - t\|\nabla V_t(y(k-1))\|^2 \\
    \geq& f(x(k)) + V_t(y(k)) + \nabla V_t(y(k))^\top (y(k-1)-y(k)) \\
    &+ \frac{t}{2} \|\nabla V_t(y(k-1))-\nabla V_t(y(k))\|^2 - t\|\nabla V_t(y(k-1))\|^2 \\
    =& \underbrace{f(x(k)) + V_t(y(k)) + \nabla V_t(y(k))^\top (x(k)-y(k))}_{=F_{y(k)}(x(k))} \\
    &+ \nabla V_t(y(k))^\top \underbrace{(y(k-1)-x(k))}_{=t\nabla V_t(y(k-1))} \\ 
    &+ \frac{t}{2} \|\nabla V_t(y(k-1))-\nabla V_t(y(k))\|^2 - t\|\nabla V_t(y(k-1))\|^2 \\
    =& F_{y(k)}(x(k)) + \frac{t}{2} \|\nabla V_t(y(k))\|^2 -\frac{t}{2} \|\nabla V_t(y(k-1))\|^2
\end{align*}
    is obtained from \eqref{eq:smooth2}. Hence, \eqref{eq:lem_acc_2} holds. $\square$
\end{pf}

With this in mind, we consider the following update rule with $\hat{x}(0) = x(0)$ and some $\{\theta_k\}\subset \mathbb{R}$:
\begin{align}\label{nesterov_theta}
    y(k) &= \hat{x}(k) - t \nabla f (\hat{x}(k)) \nonumber \\
    x(k+1)  &= T(y(k))  \nonumber \\
    \hat{x}(k+1) &=  x(k+1) + \theta_k (x(k+1)-x(k)).
\end{align}
In addition, we define $H_k:\BRnd\to\mathbb{R}$ as 
\begin{equation}\label{def:H}
    H_k = F_{y(k-1)}(x(k)) + \frac{t}{2}\|\nabla V_t(y(k-1))\|^2.
\end{equation}
with $F_y$ in \eqref{def:F_y}.
By $x(k)-y(k-1) = -t \nabla V_t(y(k-1))$, 
$H_k$ can be rewritten as $H_k = f(x(k))+V_t(y(k-1))- \frac{1}{2t} \|y(k-1)-T(y(k-1))\|^2 = f(x(k))+V_t(y(k-1)) - \frac{1}{2t} \|y(k)-x(k+1)\|^2$.

Remarkably, $H_k$ in \eqref{def:H} satisfies the following lemma.
\begin{lem}
    Consider the sequence generated by \eqref{nesterov_theta}. Then, it holds that
    \begin{equation}\label{eq:bound_H_k}
    f(x(k)) + V_t(x(k)) \leq H_k.
    \end{equation}
\end{lem}
\begin{pf}
    In light of $1/t$-smoothness of $V_t$ and $\nabla V_t(y(k-1))= -(y(k-1)-x(k))/t$, we obtain $V_t(x(k)) 
        \leq V_t(y(k-1)) 
        + \nabla V_t(y(k-1))^\top (y(k-1)-x_{k})  + \frac{1}{2t} \|y(k)-x(k+1)\|^2 
        = V_t(y(k-1)) - \frac{1}{2t} \|y(k)-x(k+1)\|^2.$
    Hence, adding $f(x(k))$ to both sides yields \eqref{eq:bound_H_k}. 
    $\square$
\end{pf}
Furthermore,
for $H_k$ in \eqref{def:H}, we prove the following inequality, which is essential to the proof of Theorem \ref{theorem:non_acc}c and \ref{theorem:acceleration}.
\begin{lem}
    For the sequence generated by \eqref{nesterov_theta} and $H_k$ defined in \eqref{def:H}, it holds that
    \begin{align}\label{eq:lem_H_k}
        H_k-H_{k+1} &\geq \frac{1}{2t} \| \hat{x}(k)- x(k+1) \|^2  \nonumber \\
        &+ \frac{1}{t} (x(k+1)-\hat{x}(k))^\top (\hat{x}(k)-x(k)).
    \end{align}
\end{lem}
\begin{pf}
    Substituting $x=x(k+1),\,y=y(k),$ and $z=x(k)$ into \eqref{eq:lem_acc_1}, we obtain
    \fontsize{9.8pt}{2.5pt}\selectfont
    \begin{align*}
        &H_{k+1} = f(x(k+1)) + V_t(y(k)) \\
        &+ \nabla V_t(y(k))^\top (x(k+1)-y(k)) + \frac{t}{2}\|\nabla V_t(y(k))\|^2 \\
        &\leq f(x(k)) + V_t(y(k)) \\
        &+ \nabla V_t(y(k))^\top (x(k)-y(k)) + \frac{t}{2} \|\nabla V_t(y(k))\|^2 \\
        &+ \frac{1}{t} (\hat{x}(k)-x(k+1))^\top (\hat{x}(k)-x(k)) - \frac{1}{2t} \|\hat{x}(k)-x(k+1)\|^2 \\
        &= F_{y(k)}(x(k)) + \frac{t}{2}\|\nabla V_t(y(k))\|^2 \\
        &+  \frac{1}{t} (\hat{x}(k)-x(k+1))^\top (\hat{x}(k)-x(k)) 
        -\frac{1}{2t} \|\hat{x}(k)-x(k+1)\|^2 \\
        &\leq F_{y(k-1)}(x(k)) + \frac{t}{2}\|\nabla V_t(y(k-1))\|^2 \\
        &+  \frac{1}{t} (\hat{x}(k)-x(k+1))^\top (\hat{x}(k)-x(k)) 
        -\frac{1}{2t} \|\hat{x}(k)-x(k+1)\|^2 \\
        &= H_k +  \frac{1}{t} (\hat{x}(k)-x(k+1))^\top (\hat{x}(k)-x(k)) \\
        &-\frac{1}{2t} \|\hat{x}(k)-x(k+1)\|^2.
    \end{align*}
    from \eqref{eq:convex_inquality}, \eqref{eq:smooth1}, and \eqref{eq:lem_acc_2}. Thus, \eqref{eq:lem_H_k} holds. $\square$
\end{pf}
\normalsize
Moreover, for the relationship between $x(k)$ and an optimal solution $x_*$, we present the following lemma.
\begin{lem}
    For $x_*\in \argmin_{x\in\CD} f(x)$, it holds that
    \begin{align}\label{eq:lem_H_k_2}
    f(x_*)+V_t(x_*) &- H_{k+1} \geq
    \frac{1}{2t} \|\hat{x}(k)-x(k+1)\|^2
    \nonumber \\
    +& \frac{1}{t} (x(k+1)-\hat{x}(k))^\top (\hat{x}(k)-x_*).
    \end{align}
\end{lem}
\begin{pf}
Recalling \eqref{nesterov_theta}, $L$-smoothness of $f$, and $1/t$-smoothness of $V_t$ for $t\in(0,1/L]$, we obtain
\fontsize{9.8pt}{2.5pt}\selectfont
\begin{align*}
    &H_{k+1} \leq f(\hat{x}(k)) - \nabla f(\hat{x}(k))^\top (\hat{x}(k)-x(k+1))\\
    & + \frac{1}{2t} \|\hat{x}_k-x(k+1)\|^2 + V_t(y(k)) - \frac{1}{2t} \|y(k)-T(y(k))\|^2 \\
    &\leq f(x_*) + \nabla f(\hat{x}(k))^\top (\hat{x}-z)- \nabla f(\hat{x}(k))^\top (\hat{x}-T(y(k)))  \\
    & + \frac{1}{2t} \|\hat{x}(k)-T(y(k))\|^2+V_t(x_*) - \frac{1}{2t} \|y(k)-T(y(k))\|^2 \\
    & + \frac{1}{t}(y(k)-T(y(k))^\top (T(y(k))-x_*+y(k)-T(y(k))) \\
    &-\frac{1}{2t} \|y(k)-T(y(k)) - (x_*-T(x_*))\|^2  \\
    &=f(x_*) + V_t(x_*) + \frac{1}{t} (\hat{x}(k)-x(k+1))^\top (\hat{x}(k)-x_*) \\
    &\quad - \frac{1}{2t} \|\hat{x}(k)-x(k+1)\|^2,
\end{align*}
\normalsize
from \eqref{eq:convex_inquality}, \eqref{eq:smooth1}, and \eqref{eq:smooth2},
where the last line is obtained because $x_*=T(x_*)$ holds for $x_*\in\CD$.
Therefore, \eqref{eq:lem_H_k_2} is obtained.
$\square$
\end{pf}

\subsection{Proof of Theorem \ref{theorem:non_acc}c}\label{subsec:proof_fixed}
In this proof, assume that $\theta_k=0$ for all $k$.
Then, $\hat{x}(k)=x(k)$ holds and the algorithm in \eqref{nesterov_theta} equals to the CPGD with $\lambda_k=t\in(0,1/L]$ for all $k\in\mathbb{N}$.

In light of \eqref{eq:lem_H_k_2} and $\hat{x}(k)=x(k)$, we obtain $2t(H_{k+1}-f^*+V_t(x_*))\leq\|x_*-x(k)\|^2$ because $2t(H_{k+1}-f^*+V(x_*)) \leq  2(x(k)-x(k+1))^\top (x(k)-x_*)
- \frac{1}{2t} \|x(k)-x(k+1)\|^2 
= \|x_*-x(k)\|^2 - \|x_*-x_{k+1}\|^2 \leq\|x_*-x(k)\|^2.$
Besides, invoking \eqref{eq:lem_H_k}, we have
\begin{equation*}
    2t(H_{k+1}-H_k )\leq \|x(k)-x(k+1)\|^2 \leq 0.
\end{equation*}
Then, following the same procedure as Theorem 3.1 in \cite{Beck2009-bq} (with $F(\bx_k) \equiv H_k$ and $F(\bx^*) \equiv f^*+V_t(x_*)$) and using \eqref{eq:bound_H_k}, we obtain \eqref{eq:fixed_bound}. $\square$

\subsection{Proof of Theorem \ref{theorem:acceleration}}\label{subsec:proof_acc}
Substituting $\theta_k=(\sigma_k-1)/\sigma_{k+1}$ in \eqref{eq:nesterov} into \eqref{nesterov_theta} yields the ACPGD in \eqref{eq:nesterov}.

Now, by \eqref{eq:lem_H_k}, \eqref{eq:lem_H_k_2}, and $\sigma_{k-1}^2 = \sigma_k(\sigma_k-1) $, following the procedure of the proof for Theorem 4.4 in \cite{Beck2009-bq} gives
\begin{align*}
    &\sigma_{k-1}^2(H_{k}-(f^*+V_t(x_*))) - \sigma_{k}^2(H_{k+1}-(f^*+V_t(x_*)))\\
    \leq& \frac{1}{2t} (\|w_{k+1}\|^2-\|w_k\|^2),
\end{align*}
where $w_k=\sigma_k (\hat{x}(k)-x_*) - (\sigma_k-1)(x(k)-x_*)$.
Thus, summing both sides over $k=1,2,\ldots$ yields
\begin{equation*}
  \sigma_k^2(  H_{k+1} - (f(x_*)+V_t(x_*))) \leq \frac{1}{2t}\|w_0\|^2 =\frac{1}{2t}\|x_0-x_*\|^2.
\end{equation*}
By $\sigma_k\geq (k+1)/2$, which can be shown by mathematical induction, we obtain
\begin{equation*}
    H_{k+1} - (f(x_*)+V_t(x_*)) \leq \frac{2\|x_0-x_*\|^2}{t(k+1)^2}. 
\end{equation*}
Therefore, the inequality \eqref{eq:nesterov_bound} follows from \eqref{eq:bound_H_k}. $\square$








\end{document}